\documentclass[10pt, a4paper]{amsart} % font to be 10
\usepackage[hmargin=20mm, vmargin=20mm, includefoot, twoside]{geometry} 
\usepackage[bookmarksopen=true, colorlinks=true, urlcolor = blue]{hyperref}
\usepackage[cp1250]{inputenc}  
\usepackage[english]{babel}
\usepackage{amssymb, verbatim, textcomp, fontenc, euscript, pstricks, pst-plot, latexsym, mathrsfs, enumerate, paralist}
\usepackage[alphabetic]{amsrefs}
\newtheorem{thm}{Theorem}
\newtheorem{step}{Step}
\newtheorem{fact}{Fact}
\theoremstyle{definition}
\newtheorem{remark}[thm]{Remark}
\newtheorem{defn}[thm]{Definition}
\begin{document}

%%%%%%%%%%%%%%%%%%%%%%%%%%%%%%%%%%%%%%%%%%%%%%%%%%%%%%%%%%%%%%%%%%%%%%%%%%%%%%%%%%%%
% Head = title, author, institute, email, grant, date, keywords, subjclass, abstract
%%%%%%%%%%%%%%%%%%%%%%%%%%%%%%%%%%%%%%%%%%%%%%%%%%%%%%%%%%%%%%%%%%%%%%%%%%%%%%%%%%%%
\title{Finitely presented subgroups of systolic groups are systolic}
\author[G.~Zadnik]{Ga\v{s}per Zadnik{$^\star$}}
\address{In\v{s}titut za matematiko, fiziko in mehaniko, Jadranska ulica 19, SI-1111 Ljubljana, Slovenija and
Institute of Mathematics of the Polish Academy of Sciences, ul. \'{S}niadeckich 8, 00-956 Warsaw, Poland}
\email{zadnik@fmf.uni-lj.si}
\thanks{$^\star$ Supported by the Slovenian Research Agency and by
the Foundation for Polish Science, project ``Cubulating groups''.}
\date{April 2013}
\keywords{Simplicial non-positive curvature, systolic complex, systolic group}
\subjclass[2010]{05E18}
\begin{abstract}
In this note we prove that every finitely presented subgroup of a systolic group is itself systolic. %We also give an example of a finitely generated subgroup of a (hyperbolic) systolic group which is not finitely presented and hence not systolic.
\end{abstract}

%%%%%%%%%%%
% make head
\maketitle

%%%%%%%
% Intro
%%%%%%%
\section{Introduction}

% To change -- only before-definition part + uncomment after-definition part
%%%%%%%%%%%%%%%%%%%%%%%%%%%%%%%%%%%%%%%%%%%%%%%%%%%%%%%%%%%%%%%%%%%%%%%%%%%%
% or copy before-theorem from FP...[fin1]
%%%%%%%%%%%%%%%%%%%%%%%%%%%%%%%%%%%%%%%%%%%%%%%%%%%%%%%%%%%%%%%%%%%%%%%%%%%%
In early eighties, Gromov deduced several properties of Riemannian manifolds of non-positive sectional curvature without using Riemannian structure, but only the property of the induced distance function, which he called \emph{CAT(0) inequality} \cite{BGS}. Gromov proved that for a cube complex equipped with the piecewise Euclidean metric, one can locally check CAT(0) condition in terms of combinatorial structure of the complex, see \cite[Theorem~II.5.20]{BH}.

In \cite{H,JS} there was introduced the following simplicial analogue of CAT(0) spaces. It is called simplicial non-positive curvature. 

\begin{defn} % flag, diagonal, m-large, m-systolic complex, systolic group %
A simplicial complex is \emph{flag} if every finite set of vertices that are pairwise connected by edges spans a simplex. A \emph{loop of length $m$} in a simplicial complex $X$ is a simplicial embedding of an $m$-cycle into $X$. An edge connecting two non-consecutive vertices of a loop is called a \emph{diagonal}. The property that every loop of length at least four and less than $m$ has a diagonal is called $m$-\emph{largeness}. Let $m\geq 6$. A simply connected $m$-large flag simplicial complex is called $m$-\emph{systolic}. We write only \emph{systolic} instead of 6-systolic. A group acting properly and cocompactly by automorphisms on a systolic complex is called \emph{systolic}.
\end{defn}

Note that this definition of systolic complex differs from the original one, but is equivalent \cite[Fact~1.2.(4) and Corollary~1.6]{JS}.

The purpose of this note is to prove the following theorem.\footnote{The author was told that the result was also proven independently in \cite{HP}.}

\begin{thm}\label{main} % The theorem
Any finitely presented subgroup of a systolic group is systolic.
\end{thm}

Theorem~\ref{main} was proven by Wise for torsion-free systolic groups, see \cite[\S5]{W}. Wise considers the quotient of the systolic complex under the group action, and his proof does not generalize to groups with torsion. Note that Theorem~\ref{main} is not true if we replace ``systolic'' with ``CAT(0)'' in the statement, see \cite[\S2.3.3]{Bra}.

%%%%%%%%%%%%%%%%%%%%%%%%%%%%%%
% Idea of the proof + notation
%%%%%%%%%%%%%%%%%%%%%%%%%%%%%%
\subsection{Notation and outline}

All the paths in any simplicial complex are taken in its one-skeleton. We use the symbol $d_X$ to denote the shortest-path distance in one-skeleton of a simplicial complex $X$, where the length of each edge is 1. Given a subcomplex $Z\leq X$, the \emph{$r$-neighborhood} of $Z$ in $X$ is defined as the simplicial span of all vertices $r$-close to $Z$, i.e.
% neighborhood
$$N_X^r(Z) = \mathrm{Span}\left\{x\in X^{ 0 }\textrm{ }|\textrm{ }d_X(x,Z^{ 0 })\leq r\right\}.$$ 
A neighborhood of a single vertex will also be called a \emph{ball} around that vertex. 
% group G and subgroup H, complex X, graph \Gamma, notation for paths
Let $G$ be a group acting properly and cocompactly by automorphisms on a systolic complex $X$. Let $H\leq G$ have a finite presentation $\langle \mathcal{S} | \mathcal{R}\rangle$ with $\mathcal{S}$ symmetric. Let $\mathcal{C}_\mathcal{S}(H)$ be the oriented Cayley graph of $H$ with respect to the generating set $\mathcal{S}$. This means that an edge connecting $h$ and $hs$ for $h\in H$ and a generator $s\in\mathcal{S}$ comes equipped with two orientations, one for $s$ and one for $s^{-1}$, except when $s^2=\mathbf{1}$, when there are two edges connecting $h$ and $hs$. Denote by $\mathcal{C}^{X}_\mathcal{S}(H)$ such subdivision of $\mathcal{C}_\mathcal{S}(H)$ that there exists a simplicial $H$-equivariant map $\phi:\mathcal{C}^{X}_\mathcal{S}(H)\to X$. Let $e_s$ denote the path between $\mathbf{1}$ and $s\in\mathcal{S}$ in $\mathcal{C}^{X}_\mathcal{S}(H)$ that comes from subdivision of the edge connecting $\mathbf{1}$ and $s$ in $\mathcal{C}_\mathcal{S}(H)$. Let $x_0=\phi(\mathbf{1})$ and $\gamma_s=\phi(e_s)$ for $s\in\mathcal{S}$. Denote by $L$ the maximum of the lengths of $\gamma_s$. Denote also $\Gamma = \phi(\mathcal{C}^{X}_\mathcal{S}(H))$. We will frequently use $s_1,\dots, s_m$ to denote generators from $\mathcal{S}$. We write $\gamma_{s_1\cdots s_m}$ for the path which is the concatenation $\gamma_{s_1} * \left(s_1\gamma_{s_2}\right) * \cdots * \left(s_1 \cdots s_{m-1}\gamma_{s_m}\right)$.

% outline
The outline of our proof is as follows. We find a neighborhood $N$ of $\Gamma$ in $X$ such that every loop in $\Gamma$ can be contracted in $N$. But new loops can appear in $N$. We thus have to consider $Y$, a disjoint union of $H$-translates of some large ball in $X$ modulo appropriate equivalence relation. Then $Y$ encodes an appropriate neighborhood of $\Gamma$ and moreover does not give rise to any new homotopically nontrivial loop. Finally, we extend $Y$ in appropriate cathegory to a maximal $H$-cocompact simply connected flag simplicial complex and prove that it is 6-large. 

%%%%%%%
% Proof
%%%%%%%
\section{Proof of the theorem}

We proceed in several steps as mentioned above. In the first step, we find a constant $R$ such that loops in $\Gamma$ are homotopically trivial in $N^{R}_X(\Gamma)$. Important properties of loops in $\Gamma$ deduced in the proof of Step~\ref{nbh} are collected in Fact~\ref{short}, since we will need them later on.

%%%%%%%%%%%%%%%%%%%%%%%%%%%%%%%%%%%%%%%%%
%Step 1 - loops from Z dies in a nbh of Z
\begin{step}\label{nbh}
There exists a constant $R<\infty$ such that every loop in $\Gamma$ is homotopically trivial in $N^{R-L}_X(\Gamma)$.
\end{step} 

\begin{proof}
After replacing a loop with its $H$-translate, it is enough to consider loops in $\Gamma$ containing a point at distance at most $L$ from $x_0$. We distinguish two main cases.

% CASE 1
\textbf{Case (1).} Loop of the form $\gamma_{s_1\cdots s_m}$ with $s_1\cdots s_mx_0=x_0$. There are three subcases.
\begin{itemize}
% 1.a
	\item[\textbf{(a)}] The word $s_1\cdots s_m$ belongs to $\mathcal{R}$. Because $\mathcal{R}$ is finite, there is a number $R_1$ such that every such loop is homotopically trivial in $N_X^{R_1}(z)$ for every vertex $z\in\gamma_{s_1\cdots s_m}$. 
% 1.b	
	\item[\textbf{(b)}] The word $s_1\cdots s_m = \mathbf{1} $ but it does not belong to $\mathcal{R}$. Then $s_1\cdots s_m$ is a concatenation of conjugates of relators from $\mathcal{R}$, but each such conjugate is homotopically trivial in $N_X^{R_1}(\Gamma)$ hence the whole loop is homotopically trivial in $N_X^{R_1}(\Gamma)$.
% 1.c	
	\item[\textbf{(c)}] The point $x_0$ is fixed by $s_1\cdots s_m$ but $s_1\cdots s_m\neq\mathbf{1}$. Without loss of generality, we can assume that $s_1\cdots s_m$ is the shortest representative of the corresponding group element since all the other representatives differ from the shortest by concatenation with words considered in Subcases~(1.a,~1.b). By properness of $G$-action and hence of $H$-action, the number of elements $h\in H$ fixing $x_0$ is finite. Hence we can choose a constant $R_1'$ such that $\gamma_{s_1\cdots s_m}$ is homotopically trivial in $N_X^{R'_1}(z)$ for every vertex $z\in\gamma_{s_1\cdots s_m}$. 
\end{itemize}

% CASE 2	
\textbf{Case (2).} Loop coming from a path $\gamma_{s_1\cdots s_m}$ with self-intersection $x\notin Hx_0$. Without loss of generality we can assume that $x=\gamma_{s_1}\cap (s_{1}\cdots s_{m-1} \gamma_{s_m})$. Figure \ref{skica1} shows such configuration. Observe that in this case, $d_X(x_0,s_1\cdots s_m x_0)\leq 2L$. By properness of the $H$-action, there is an upper bound $N$ such that $s_1\cdots s_m = p_1\cdots p_k$, where the number $k$ of terms $p_i\in\mathcal{S}$ is at most $N$.  Since 
$$s_1\cdots s_mp_k^{-1}\cdots p_1^{-1} = \mathbf{1},$$ 
the big loop $\gamma_{s_1\cdots s_mp_k^{-1}\cdots p_1^{-1} }$ is a loop from Cases~(1.a,~1.b), hence it can be contracted in $N_X^{R_1}(\Gamma)$. Thus to contract the original loop it suffices to contract 
%%%%%%%%%%%%%%%%%%%%%%%%%%%%%%%%%%%%%%%%%%%%%%%%%
\begin{equation}\label{loop2} % loops of type (2)
\gamma_{p_1\cdots p_k} * \gamma, \tag{$\diamondsuit$}
\end{equation} 
where $\gamma$ is a path in $\Gamma$ from $s_1\cdots s_mx_0$ to the intersection $x$ of $\gamma_{s_1} $ and $s_1\cdots s_{m-1}\gamma_{s_m}$ concatenated with a path from $x$ to $x_0$. But the total length of the loop (\ref{loop2}) is at most $(N+2)L$. Let $R_2$ be a number such that any loop $\gamma_{p_1\cdots p_k} * \gamma$ of type (\ref{loop2}) can be contracted in $N_X^{R_2}(z)$ for any vertex $z\in\gamma_{p_1\cdots p_k} * \gamma$.

% SKETCH for case 2
\begin{figure}[h!]
	\centering
	{\psset{unit=0.6cm}
	\begin{pspicture}(-1,-1.5)(10,3)	
% path x_0 - p_1x_0 - ... - p_k...p_1x_0
	{\red \pscurve[linecolor=red](0,0)(-1,-.5)(-2,1)(-1,2)(1,2)
	\psdots[linecolor=red](-1,-.5)(-1,2)
	\uput[dl](-1,-.4){\tiny{$p_1x_0$}}
	\uput[ul](-1,1.9){\tiny{$p_1\cdots p_{k-1}x_0$}}
	\uput[l](-2,1){\tiny{$\gamma_{p_1\cdots p_k}$}}}
% some points from H-orbit of 0-skeleton
	\psdots(0,0)(5,1)(4,-1)(1,2)
	\psline(0,0)(2,1)(4,0)(5,1)
	\pscurve(5,1)(7,2)(10,2)(12,1)(11,-1)(7,-2)(4,-1)
	\psline(4,-1)(3,-1)(2,1)(1,2)
% path p
	{\green \uput[u](1,.5){\tiny{$\gamma$}}}
	\psline[linecolor=green,linestyle=dotted, linewidth=.15](1,2)(2,1)(0,0)
% denotation	
	\uput[ur](2,1){\tiny{$x$}}
	\uput[d](0,0){\tiny{$x_0$}}
	\uput[r](5,1){\tiny{$s_1x_0$}}
	\uput[u](1,1.9){\tiny{$s_1\cdots s_mx_0$}}
	\uput[u](4.2,-1.1){\tiny{$s_1\cdots s_{m-1}x_0$}}
	\uput[r](12,1){\tiny{$\gamma_{s_1\cdots s_m}$}}
	
	\end{pspicture}}
	\caption{Loop from Case~(2); possibly $x=x_0$ or $x=s_1\cdots s_mx_0$. If both equalities hold, this example is covered in Case~(1).}	
	\label{skica1}
\end{figure}
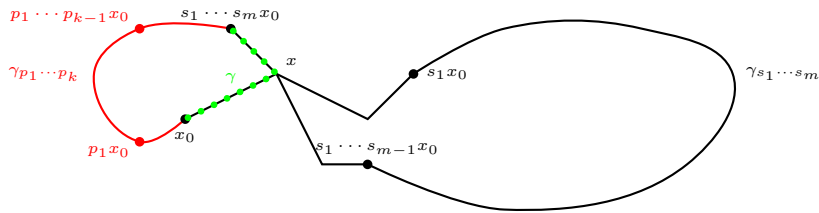

Alltogether, we find a constant $R' = \max\{ R_1,R'_1,R_2\}$ such that every loop in $\Gamma$ is homotopic in $N^{R'}_X(\Gamma)$ to a trivial loop. Hence $R = R'+L$ works. 
\end{proof}

We will call loops from Cases~(1.a,~1.c) and loops (\ref{loop2}) form Case~(2) \emph{short loops}. From the proof of Step~\ref{nbh} we deduce the following.

\begin{fact}\label{short}
Every loop in $\Gamma$ is a concatenation of conjugates of short loops. Let $\gamma$ be a short loop in $\Gamma$. Then for every vertex $z$ on $\gamma$, the loop $\gamma$ is fully contained in the ball $N^{R-L}_X(z)$. \hfill\qed
\end{fact}

%%%%%%%%%%%%%%%%%%%%%%%%%%%%%
% Step 1+ - construction of Y
We are now ready to define $Y$. For every $h\in H$, denote $B_h^{ 0 }=N^{R}_X(hx_0)^{ 0 }$ and denote the copy of the vertex $v\in X$ in $B_h^{ 0 }$ by $v^h$. Let $\sim$ be an equivalence relation on $\coprod_{h\in H}B_h^{ 0 }$ generated by $v^h\sim u^g$ if and only if $v=u$ in $X$ and $g^{-1}h\in\mathcal{S}$. Let 
$$Y^{ 0 }=\left(\coprod_{h\in H} B_h^{ 0 }\right)\Big/\sim.$$ 
Note that $B_h^{ 0 }$ injects into $Y^{ 0 }$. For $y\in Y^{ 0 }$ we write $\bar{y}$ for the vertex of $X$ such that $y=\bar{y}^h$ for some $h\in H$.
Next, we define $Y^{ 1 }$. We connect two vertices $y,z\in Y^{ 0 }$ by an edge if there exist representatives $\bar{y}^h$, $\bar{z}^g$ for $y$ and $z$ with $h=g$ and $\bar{y},\bar{z}$ adjacent. Let $Y$ be the flag completion of $Y^{ 1 }$ and let $B_h$ be the simplicial span of $B_h^{ 0 }\leq Y$. We define a natural action of $H$ on $Y$, which is induced from $H$-action on $X$. Note that $hB_\mathbf{1} = B_h$, hence $H$-action on $Y$ is proper and cocompact. 

% note : proper H-map Y -> X
Observe that by the construction of $Y$, there exists a proper $H$-equivariant map $f:Y\to X$. It is defined by $f(y)=\bar y$ for $y\in Y^{ 0 }$ and extends to higher-dimensional simplices simplicially. Let us define local sections
%%%%%%%%%%%%%%%%%%%%%%%%%%%%%%%%%%%%%%%%%%%%%%%%%%%%%%%%%%%%%%%%%%%%%%%%%%%
\begin{equation}\label{i-h-maps} % the maps i_h from balls in X to B_h in Y
i_h:N_X^R(hx_0)\to B_h\tag{$\heartsuit$}
\end{equation} 
by $i_h(u) = u^h$ for every $u\in N_X^R(hx_0)^{ 0 }$ and every $h\in H$. By definition of $Y$, each map $i_h$ is bijective on zero-skeleta. Furthermore, vertices $u^h,v^h\in B_h^{ 0 }$ are adjacent if and only if $u,v\in N_X^R(hx_0)^{ 0 }$ are adjacent. Hence $i_h$ is well defined isomorphism between one-skeleta of $N_X^R(hx_0)$ and $B_h$. Since a flag complex is determined by its one-skeleton, $i_h$ is an isomorphism. Note that $B_h$ might not be a ball in $Y$. 

Observe that for any $h\in H$ and $s\in\mathcal{S}$ the two maps $i_h$ and $i_{hs}$ agree on $N_X^R(hx_0)\cap N_X^R(hsx_0)$ because they agree on the zero-skeleton of that intersection. By that fact, there is a natural map $\varphi:\mathcal{C}^{X}_\mathcal{S}(H)\to Y$, sending edge path $he_s$ to $i_h(h\gamma_s)$. We can as well describe the map $\varphi$ in terms of the $H$-action on $Y$, but the previous definition is more useful for us. Next step ensures that $\varphi$ has good properties.

%%%%%%%%%%%%%%%%%%%%%%%%%%
% Step 2 - map \Gamma -> Y
\begin{step}\label{embedding}
The map $\varphi$ factors through $\phi:\mathcal{C}^{X}_\mathcal{S}(H)\to\Gamma$.
\end{step}

\begin{proof}
We have to check that if two points of $\mathcal{C}^{X}_\mathcal{S}(H)$ are identified under $\phi$, they are also identified under $\varphi$. To see this, observe that if two paths $\gamma_s$ and $h\gamma_{s'}$ in $\Gamma$ cross, where $s,s'\in\mathcal{S}$ and $h\in H$, then there is a sequence of generators $p_1,\dots, p_k\in\mathcal{S}$ such that $p_1\cdots p_k$ is the shortest word representing $h$. In particular $p_1\cdots p_ke_{s'} = he_{s'}$. %Denote $p_0=\mathbf{1}$. 
It follows from Fact~\ref{short} that the ball $N_X^{R-L}(p_1\cdots p_lx_0)$ contains the whole path $\gamma_{p_1\cdots p_k}$ for all $l=0,1,\dots,k$, where the empty word represents $\mathbf{1}$. Thus if we write $\overline{\gamma_s}$ for $\gamma_s$ with opposite orientation, we have that $N_X^{R}(p_1\cdots p_lx_0)$ contains the path $\gamma=(h\gamma_{s'})*\gamma_{p_1\cdots p_k}*\overline{\gamma_{s}}$ for all $l=0,1,\dots,k$. Hence $B_{p_1\cdots p_l}$ contains $i_{p_1\cdots p_l}(\gamma)$ for all $l=0,1,\dots,k$. Since two consecutive maps $i_{p_1\cdots p_{l-1}}$ and $i_{p_1\cdots p_l}$ agree on the intersection of their domains, the path $i_\mathbf{1}(\gamma)=i_h(\gamma)$ is contained in the intersection $\bigcap_{l=0}^kB_{p_1\cdots p_l}$. Thus the point on $e_s$ is identified with the appropriate point on $he_{s'}$ under $\varphi$. This finishes the proof. 
\end{proof}

% section \Gamma -> Y
By Step~\ref{embedding}, there exists a lift $f_\Gamma:\Gamma\to Y$ of the map $f:Y\to X$. Obviously $f_\Gamma$ agrees with $i_h$ on $\Gamma\cap N_X^R(hx_0)$. From now on, we identify $\Gamma$ with its $f_\Gamma$-image.

%%%%%%%%%%%%%%%%%%%%%%%%%%%
% Step 3 - Y is 1-connected
\begin{step}\label{1connected}
The complex $Y$ is simply connected.
\end{step}

As mentioned in the outline, we first prove that $Y$ encodes an appropriate neighborhood of $\Gamma$ such that loops in $\Gamma$ are homotopically trivial in $Y$. Then we exhibit a homotopy from any loop in $Y$ to a loop in $\Gamma$.

\begin{proof} % loops in \Gamma are trivial
Take any loop $\gamma$ in $\Gamma\subseteq Y$. By Fact~\ref{short}, it is a concatenation of short loops. In the same way as in the proof of Step~\ref{embedding} one can show that each short loop $\gamma'$ is fully contained in $B_h$ for each $h\in H$ such that $d_Y(hx_0,\gamma')\leq L$. Pick such $h\in H$. We know that $B_h$ is isomorphic to $N_X^R(hx_0)$ via the map $i_h$ from (\ref{i-h-maps}). Invoking Fact~\ref{short} once again, we see that the loop $\gamma'$ is homotopically trivial in $B_h$, hence $\gamma$ is homotopically trivial in $Y$.

% Homotopy from a loop in Y to a loop in \Gamma
Finally we need to show that every loop in $Y$ is homotopic to a loop in $\Gamma$. Let $\beta:S^1\to Y^{ 1 }$ be a loop in the one-skeleton of $Y$. We identify $S^1$ with $I/\partial I$, where $I=[0,1]$ is an interval. Let $0\leq t_0<t_1<\dots<t_n<1$ be cyclically ordered points on $S^1$ and $h_0,h_1,\dots,h_n$ elements of $H$ such that $\beta(t_i)\in Y^{ 0 }$ and $\beta([t_{i},t_{i+1}])\in B_{h_i}$ for all $i=0,1,\dots, n$, where indices are taken modulo $n+1$. Since $B_{h_{i-1}}$ and $B_{h_{i}}$ both contain $\beta(t_i)$, there is a sequence of generators $s^i_1,\dots, s^i_{n(i)}\in\mathcal{S}$ such that $h_i = h_{i-1}s^i_1\cdots s^i_{n(i)}$ and $\beta(t_i)\in B_{h_{i-1}s^i_1\cdots s^i_l}$ for all $l=0,1,\dots, n(i)$. Recall that empty word stands for $\mathbf{1}$. This means that there exist geodesics 
% paths \beta to define a homotopy
$$\beta_l^i:\left(I,\partial I\right)\to\left(B_{h_{i-1}s_1\cdots s_l},\{\beta(t_i),h_{i-1}s_1\cdots s_lx_0\}\right)\textrm{ for all } i=0,1,\dots,n\textrm{ and }l=0,1,\dots,s(i).$$
Recall that $h_{i-1}s_1\cdots s_lx_0$ belongs also to $B_{h_{i-1}s_1\cdots s_{l-1}}$. Since balls in systolic complexes are geodesically convex \cite[\S7]{JS}, the image of $\beta_l^{i}$ is contained in $B_{h_{i-1}s_1\cdots s_{l-1}}$. Next, we can find some $\varepsilon>0$ such that $t_i+n(i)\varepsilon<t_{i+1}$ for all $i$. After precomposing $\beta$ with a map $S^1\to S^1$, homotopic to the identity, we can assume that $\beta$ is constant on $[t_i,t_i+n(i)\varepsilon]$ for all $i=0,1,\dots,n$.  Hence
\begin{itemize}
	\item for all $i=0,1,\dots,n$ and $l=1,2,\dots,n(i)$ there is $H_l^i:[t_i+(l-1)\varepsilon,t_i+l\varepsilon]\times I \to B_{h_{i-1}s_1\cdots s_{l-1}}\cap B_{h_{i-1}s_1\cdots s_l}$ with $H_l^i(t_i+(l-1)\varepsilon,t) = \beta_{l-1}^{i}(t)$ and $H_l^i(t_i+l\varepsilon,t) = \beta_l^{i}(t)$, where $H_l^i(-,0)$ is contant path $\beta(t_i)$ and $H_l^i(-,1)$ is a path $h_{i-1}s_1\cdots s_{l-1}\gamma_{s_l}\subseteq\Gamma$;
	\item for all $i=0,1,\dots,n$ there is $H^i:[t_i+n(i)\varepsilon,t_{i+1}]\times I \to B_{h_i}$ with $H^i(t_i+n(i)\varepsilon,t) = \beta_{n(i)}^i(t)$ and $H^i(t_{i+1},t) = \beta_{0}^{i+1}(t)$, where $H(-,0) = \beta|_{[t_i+n(i)\varepsilon,t_{i+1}]}$ and  $H(-,1)$ is a constant path $h_ix_0$.
\end{itemize}
Since the homotopies from above agree on the intersections of their domains, they glue together to a homotopy $H:S^1\times I\to Y$ with $H(-,0) = \beta$ and $H(-,1)$ a loop in $\Gamma$.
\end{proof}

% Last step
In the following step, using $f$ we extend $Y$ to a systolic complex $\overline{Y}$, on which $H$ still acts properly and cocompactly and is thus a systolic group. 

% Family of f-extensions of Y; maximal element is the space we want
We say that a pair $(W,f_W)$ is an \emph{$f$-extension of $Y$} if the following holds. The complex $W$ is simply connected flag simplicial complexes containing $Y$ such that $Y^{ 0 } = W^{ 0 }$ and that the $H$-action on $Y$ extends to an $H$-action on $W$. Furthermore, the map $f_W: W\to X$ is a simplicial $H$-equivariant map which extends $f$. Note that $f_W$ maps an edge of $W$ either to an edge or to a vertex of $X$.

Let $\mathcal{F}$ be the family of all $f$-extensions of $Y$. Observe that $\mathcal{F}$ is equipped with a natural partial order $\leq$, where $(W_1,f_{W_1})\leq (W_2,f_{W_2})$ if there exists an $H$-equivariant embedding $i : W_1\to W_2$ fixing $Y$ such that $f_{W_2}\circ i = f_{W_1}$. The family $\mathcal{F}$ is nonempty since it contains $(Y,f)$ by Step~\ref{1connected}. Let $\left(W_\lambda,f_{W_\lambda}\right)_{\lambda\in\Lambda}$ be an increasing chain in $\mathcal{F}$. Then the union $\left(\bigcup_{\lambda}W_\lambda,\bigcup_\lambda f_{W_\lambda}\right)$ is also in $\mathcal{F}$, so it is an upper bound for the chain $(W_\lambda,f_{W_\lambda})_{\lambda\in\Lambda}$. By Kuratowski--Zorn Lemma, there exists a maximal element $(\overline{Y},\overline{f})\in \mathcal{F}$. 

%%%%%%%%%%%%%%%%%%%%%%%%%%%%%%%%%%%
% Step 4 - \overline{Y} is systolic
\begin{step}
For a maximal element $(\overline{Y},\overline{f})\in \mathcal{F}$, the simplicial complex $\overline{Y}$ is a systolic complex, equipped with a proper and cocompact $H$-action.
\end{step}

% \overline{Y} is proper and H-cocompact
\begin{proof} We claim that the valence in $\overline{Y}^{ 1 }$ of each $y\in \overline{Y}^{ 0 }$ is bounded from above. Recall that $\overline{f}$ and $f$ agree on $\overline{Y}^{ 0 } = Y^{ 0 }$. Let $N_y\subseteq  Y^{ 0 }$ denote the set of all vertices adjacent to $y$ in $\overline{Y}$. For every $y'\in N_y$ either $f(y')=f(y)$ or $f(y')$ is adjacent to $f(y)$. This means that $f(N_y)$ is contained in $N_X^1(f(y))$.
In other words, we have $N_y\subseteq f^{-1}(N_X^1(f(y)))$. Because $X$ is proper, the ball $N_X^1(f(y))$ is compact. But $f$ is a proper map, hence the set $f^{-1}(N_X^1(f(y)))$ is compact and the claim is proven. In particular, $\overline{Y}$ is a proper simplicial complex. Because the vertex set of $Y$ and $\overline{Y}$ coincide, the action of $H$ on $\overline{Y}$ is proper and cocompact. By definition, $\overline{Y}$ is flag and simply connected. It remains to prove 6-largeness.

% \overline{Y} is 6-large
Suppose for contradiction that there is some loop $\alpha$ of length four or five in $\overline{Y}$ without diagonals. If $\overline{f}$ maps $\alpha$ bijectively to $\alpha'=\overline{f}(\alpha)\subseteq X$, then there exist two non-consecutive vertices $u'$ and $v'$ of $\alpha'$ connected by a diagonal because $X$ is systolic. Let $u$ and $v$ be the vertices of $\alpha$ mapped to $u'$ and $v'$ by $\overline{f}$. For every $h\in H$, we add an edge in $Y$ between $hu$ and $hv$ and extend $\overline{f}$ to the new edges naturally. Let us remind the reader that if for $n$ different $h_1,\dots, h_n\in H$ all the sets $\{h_iu,h_iv\}$ coincide for $i=1,2,\dots,n$, we only add one edge between $h_1u$ and $h_1v$ instead of $n$. This remark will be applied two more times without mentioning it.

If $\overline{f}$ is not bijective on $\alpha$, there must be two vertices $u$, $v$ of $\alpha$, which are mapped by $\overline{f}$ to the same vertex. If they are non-consecutive in $\alpha$, we add edges between $hu$ and $hv$ for every $h\in H$ and extend $f$ such that it maps any new edge to the common image of its endpoints. If $u$ and $v$ are consecutive, let $w\neq u$ be the other neighbor of $v$ in $\alpha$. Then we add edges between $hu$ and $hw$ for all $h\in H$. Note that since $\overline{f}(u) = \overline{f}(v)$, the point $\overline{f}(w)$ is either adjacent to or coincide with $\overline{f}(u)$ and hence we can extend $\overline{f}$ to the newly added edges.

In all cases, we added an $H$-orbit of an edge to $\overline{Y}$. After a flag completion, we obtain a flag simplicial complex $\hat{Y}$ on the set of vertices $Y^{ 0 }$, properly containing $\overline{Y}$, together with a map $\hat{f}$, extending $\overline{f}$, and equipped with an $H$-action, extending $H$-action on $\overline{Y}$. The complex $\hat{Y}$ is also simply connected. Indeed, every edge $e$ in $\hat{Y}^{ 1 }-\overline{Y}^{ 1 }$ is a diagonal of a loop $\alpha$ in $\overline{Y}^{ 1 }$ of length less than six. This means that $e$ together with two consecutive edges of $\alpha$ form a triangle, which is filled after flag completion. Hence the path $e$ is homotopic relative to its endpoints to a path of length two in $\overline{Y}$. In other words, any loop in $\hat{Y}$ is homotopic to a loop in $\overline{Y}$ and the later is simply connected since it belongs to $\mathcal{F}$. Hence $(\overline{Y},\overline{f})\lneq (\hat{Y},\hat{f})\in\mathcal{F}$, which contradicts the maximality of $(\overline{Y},\overline{f})$.
\end{proof}
	
%%%%%%%%%%%%%%%%%%%%%%%%%%%%%%%%%%%%%%%%%%%%%%%%%%%%%%%%%%%%%%%%%%
% Can not omit finitely presented assumption to finitely generated
%%%%%%%%%%%%%%%%%%%%%%%%%%%%%%%%%%%%%%%%%%%%%%%%%%%%%%%%%%%%%%%%%%
\begin{remark}
A first counterexample to Theorem~\ref{main} assuming only finitely generated instead of finitely presented subgroup is due to Stallings. Denote a free group of rank two generated by $x$ and $y$ by $\langle x,y\rangle$. In \cite{S} (see also \cite[\S2.4.2]{Bra}) Stallings proved that the kernel $K$ of the homomorphism 
$$\tau:\langle a,b\rangle\times \langle x,y\rangle\to\mathbf{Z},\quad \tau(a) = \tau(b) = \tau(x) = \tau(y)=1,$$ is finitely generated, but not finitely presentable. Hence $K$ cannot be systolic. On the other hand, the direct product of two free groups of rank two is systolic, see \cite{EP}.

Even in the case where $G$ is hyperbolic, one cannot hope for a generalization of the theorem above. By the Rips Construction \cite{R}, for any finitely presented group $Q$ and arbitrary $\lambda>0$, there exists a finitely presented $C'(\lambda)$ small cancellation group $G$ and a short exact sequence $\{\mathbf{1}\}\to N\to G\to Q\to \{\mathbf{1}\}$, where $N$ is finitely generated normal subgroup of $G$. Due to \cite{B}, the group $N$ is finitely presentable if and only if $Q$ is finite. Hence, if we choose $Q=\mathbf{Z}$ and $\lambda = \frac{1}{6}$, then the Rips Construction gives a finitely presented $C'(1/6)$ small cancellation group $G$ which is hyperbolic \cite{G} and C(7) \cite[Chapter~V,~\S2]{LS}. By \cite{W}, $C(7)$ group $G$ is 7-systolic. But it has a finitely generated not finitely presentable subgroup $N$, hence a finitely generated non-systolic subgroup. In particular, systolic and even 7-systolic groups are not coherent in general.
\end{remark}

%%%%%%%%%%%%%
% T H A N K S
%%%%%%%%%%%%%
\section*{Acknowledgement}
I would like to express my gratitude to the Institute of Mathematics of the Polish Academy of Sciences for hospitality during my work on this project, and especially to Piotr Przytycki for suggesting me this problem and for his indispensable advices and comments on previous drafts. I would also like to thank Tomasz Elsner for useful discussion on the problem. 

%%%%%%%%%%%%%%%%%%%%%%%%%%%%%%%%%
% R  E  F  E  R  E  N  C  E  S  %
%%%%%%%%%%%%%%%%%%%%%%%%%%%%%%%%%

\end{document}